\documentclass[10pt]{amsart}
\usepackage{tipa}
\usepackage{amssymb}
\usepackage{stmaryrd}
\usepackage{amsmath,amsfonts,amsthm,mathrsfs,bbm,array,subfigure}

\newtheorem{theorem}{Theorem}[section]

\theoremstyle{definition}

\newtheorem{question}[theorem]{Question}

\theoremstyle{remark}

\numberwithin{equation}{section}

\begin{document}

\title[On the Diophantine equation $f(x)f(y)=f(z)^n$]{On the Diophantine equation $f(x)f(y)=f(z)^n$ concerning Laurent polynomials}

\author{Yong Zhang}
\address{School of Mathematics and Statistics, Changsha University of Science and Technology,
Changsha 410114, People's Republic of China}

\address{School of Mathematical Sciences, Zhejiang University, Hangzhou 310027, People's Republic of China}
 \email{zhangyongzju$@$163.com}

\thanks{This research was supported by the National Natural Science Foundation of China (Grant No.~11501052).}

%    General info

\date{}


\begin{abstract}
By the theory of elliptic curves, we investigate the nontrivial
rational parametric solutions of the Diophantine equation
$f(x)f(y)=f(z)^n$, where $n=1,2$ and $f(X)$ are some simple Laurent
polynomials.

\end{abstract}

\maketitle

\section{Introduction}
Let $f(X)\in \mathbb{Q}[X]$ be a polynomial without multiple roots
and $deg(f)\geq2$. There are many authors considered the integer
solutions of the Diophantine equation
\begin{equation}\label{Eq11}
f(x)f(y)=f(z)
\end{equation} for different polynomials $f(X)$, such as L. Euler \cite{Dickson2} ($f(X)=X(X+1)/2$), C. Ko
\cite{Ko} ($f(X)=X^X$, P. Erd\H{o}s proposed), S. Katayama
\cite{Katayama1} and A. Baragar \cite{Baragar} ($f(X)=X(X+1)$), K.
Kashihara \cite{Kashihara} ($f(X)=X^2-1$), Y. Bugeaud \cite{Bugeaud}
and M.A. Bennett \cite{Bennett} ($f(X)=X^k-1,k>3$), and Y. Zhang and
T. Cai \cite{Zhang-Cai1,Zhang-Cai2} ($z$ is a square, $deg(f)=2,3$).
In 1992 and 1994, A. Schinzel and U. Zannier
\cite{Schinzel-Zannier1,Schinzel-Zannier2} investigated the number
of integer solutions of Eq. (\ref{Eq11}) for monic quadratic
polynomials with integer coefficients and gave many important
results about it. In 2015, Y. Zhang and T. Cai \cite{Zhang-Cai3}
showed that Eq. (\ref{Eq11}) has infinitely many nontrivial positive
integer solutions for $f(X)=X(X+d)$ with $d\geq3$, infinitely many
nontrivial positive integer solutions for $f(X)=(X-1)X(X+1)$, and a
rational parametric solution for $f(X)=X(X-1)(X-2)(X-3)$.

Another interesting Diophantine equation is
\begin{equation}\label{Eq12}
f(x)f(y)=f(z)^2.
\end{equation}
For the related information, we can refer to
\cite{Schinzel-Sierpinski,Szymiczek,Katayama2,Bennett,Ulas2,Guy}. In
2007, M. Ulas \cite{Ulas1} proved that if $f(X)=X^2+k,$ where $k\in
\mathbb{Z}$, then Eq. (\ref{Eq12}) has infinitely many rational
parametric solutions; if~$f(X)=X(X^2+X+t),$ where $t\in \mathbb{Q}$,
then Eq. (\ref{Eq12}) has infinitely many rational solutions for all
but finitely many $t$. In 2015, we \cite{Zhang} gave the conditions for $f(X)=AX^2+BX+C$ such that Eq.
(\ref{Eq12}) has infinitely many
nontrivial integer solutions and proved that it has a
rational parametric solution for infinitely many irreducible cubic
polynomials wihch solved
a problem of \cite{Ulas1}.

In this paper we consider the nontrivial rational parametric
solutions of Eqs. (\ref{Eq11}) and (\ref{Eq12}) for some simple
Laurent polynomials. The nontrivial solutions $(x,y,z)$ of Eq.
(\ref{Eq11}) means that $x\neq z,y\neq z$ and $f(z)\neq0$ and the
nontrivial solutions $(x,y,z)$ of Eq. (\ref{Eq12}) means that
$f(x)\neq f(y)$ and $f(z)\neq0$.

Recall that a Laurent polynomial with coefficients in a field
$\mathbb{F}$ is an expression of the form
\[f(X)=\sum_{k}a_kX^k,a_k\in \mathbb{F},\]
where $X$ is a formal variable, the summation index $k$ is an
integer (not necessarily positive) and only finitely many
coefficients $a_k$ are non-zero. Here we mainly care about the
simple Laurent polynomials
\[f(X)=AX+B+\frac{C}{X},\]
where $A,B,C\in \mathbb{Z}$ and $C\neq0$.

For Eq. (\ref{Eq11}), by the theory of elliptic curves, we have

\begin{theorem} For $f(X)=AX+B+\frac{C}{X}$ with $ABC\neq0$, Eq. (\ref{Eq11}) has infinitely many nontrivial rational parametric
solutions $(x,y,z)$.
\end{theorem}

\begin{theorem} For $f(X)=AX+\frac{C}{X}$ with $AC\neq0$, if there exists a nonzero rational number $T$ such that the elliptic curve
\[E_2:~Y^2=(X-2ACT^2)(X+2ACT^2)(X-2AC(2A^2T^4+(4AC-1)T^2+2C^2))\]
has positive rank, then Eq. (\ref{Eq11}) has infinitely many
nontrivial rational solutions $(x,y,z)$.
\end{theorem}

\textbf{Remark 1.1.} When $A=0,$ $f(X)=B+\frac{C}{X}$, Eq.
(\ref{Eq11}) becomes\[\frac{(B^2-B)xyz-Cxy+BCxz+BCyz+C^2z}{xyz}=0.\]

If $B=1,$ we can find infinitely many integer solutions of
$xy-xz-yz-Cz=0$ for any nonzero integers $C$, such as
\[(x,y,z)=(k(k+1+C),k+1,k),(k(kl-C-l),kl-C,kl-C-l),\] where $k,l$ are
integer parameters. When $f(X)=1+\frac{1}{X},$ it is interesting to
note that A. Padoa \cite[p.688]{Dickson2} studied the Diophantine
equation
\begin{equation}
\bigg(1+\frac{1}{x}\bigg)\bigg(1+\frac{1}{y}\bigg)=1+\frac{1}{z},
\end{equation}
which is equivalent with
\[(x-z)(y-z)=z(z+1).\]
If $z$ is given, then we can obtain all couples $x,y$ by finding all
pairs of positive integers whose products is $z(z+1)$, and adding
$z$ to each other. J.E.A. Steggall \cite[p.688]{Dickson2} found the
positive integer solutions of Eq. (1.3), by noting that $xy$ must be
divisible by $x+y+1=a,$ and hence $a|x(x+1).$ Hence for any integer
$x$, determine a factor $a>x+1$ of $x(x+1)$; then $y=a-x-1,$ while
$z=x-b,$ where $b=x(x+1)/a$.

If $B\neq0,1$, it is easy to get the rational parametric solutions
of $(B^2-B)xyz-Cxy+BCxz+BCyz+C^2z=0$, such as
\[(x,y,z)=\bigg(u,v,\frac{Cuv}{(B^2-B)uv+BCu+BCv+C^2}\bigg),\]where $u,v$ are rational parameters.
But it seems difficult to get its infinitely many integer solutions.

For Eq. (\ref{Eq12}), by the similar methods, we get

\begin{theorem} For $f(X)=AX+B+\frac{C}{X}$ with $ABC\neq0$, Eq. (\ref{Eq12})
has infinitely many nontrivial rational solutions $(x,y,z)$.
\end{theorem}

\begin{theorem} For $f(X)=AX+\frac{C}{X}$, if there exist nonzero integers $A,C$ and a nonzero rational number $z$ such that the quartic elliptic curve
\[C_4:~v^2=-4A^3Cz^4y^4+(A^4z^8+4A^3Cz^6-2A^2C^2z^4+4AC^3z^2+C^4)y^2-4AC^3z^4\]
has infinitely many rational points $(y,v)$, then Eq. (\ref{Eq12})
has infinitely many nontrivial rational solutions $(x,y,z)$.
\end{theorem}

\textbf{Remark 1.2.} When $A=0,$ $f(X)=B+\frac{C}{X}$, Eq.
(\ref{Eq12}) reduces to
\[\frac{C((2Byz-Bz^2+Cy)x-Byz^2-Cz^2)}{xyz^2}=0.\] If $BC\neq0$, a rational parametric solution of
$(2Byz-Bz^2+Cy)x-Byz^2-Cz^2=0$ is
\[(x,y,z)=\bigg(\frac{v^2(Bu+C)}{2Buv-Bv^2+Cu},u,v\bigg),\]where $u,v$
are rational parameters. To get its infinitely many integer
solutions is a difficult problem.

\section{The proofs of Theorems}

\begin{proof}[\textbf{Proof of Theorem 1.1.}]
For $f(X)=AX+B+\frac{C}{X}$, let \[x=T,y=u.\] Then Eq. (\ref{Eq11})
becomes
\[\begin{split}
f(u)=&Az(AT^2+BT+C)u^2+(-ATz^2+(AT^2+(B-1)T+C)Bz-CT)u\\
&+Cz(AT^2+BT+C)=0.\end{split}\]
Consider the above equation as a
quadratic equation of $u$, if it has rational solutions, then the
discriminant
\[\begin{split}
&\Delta(z)=a_1z^4+a_2z^3+a_3z^2+a_4z+a_5
\end{split}\]
should be a square, where\[\begin{split}
a_1=&A^2T^2,\\
a_2=&-2ABT(AT^2+(B-1)T+C),\\
a_3=&-A^2(4AC-B^2)T^4-2AB(4AC-B^2+B)T^3+(B^4-8A^2C^2-2AB^2C\\
&-2B^3+2AC+B^2)T^2-2BC(4AC-B^2+B)T-C^2(4AC-B^2),\\
a_4=&-2BCT(AT^2+(B-1)T+C),\\
a_5=&C^2T^2.
\end{split}\]Let us consider the curve
$C_1:~v^2=\Delta(z).$

1) If $4AC-B^2\neq0,$ the discriminant of $\Delta(z)$ is non-zero as
an element of $\mathbb{Q}(T)$, then $C_1$ is smooth. By the method described in \cite[p.77]{Mordell} (or
see \cite[p.476, Proposition 7.2.1]{Cohen1}), $C_1$ is birationally
equivalent with the elliptic curve
\[E_1:Y^2=(X+2ACT^2)(X^2+a_3X+a_6),\]
where $a_3$ is the coefficient of $z^2$ in $\Delta(z)$ and
\[\begin{split}
a_6=&2ACT^2(A^2(4AC+B^2)T^4+2AB(4AC+B^2-B)T^3+(8A^2C^2+6AB^2C\\
&+B^4-2B^3-2AC+B^2)T^2+2BC(4AC+B^2-B)T+C^2(4AC+B^2)).
\end{split}\]
Because the map $\varphi_1: C_1\rightarrow E_1$ is quite
complicated, we do not present the explicit equations for the
coordinates of it.

It is easy to see that $E_1$ contains the point
\[P=(2ACT^2,4ABC(AT^2+(B-1)T+C)T^2).\]
By the group law, we have
\[\begin{split}
[2]P=\bigg(&2AC(2A^3CT^4+4A^2C(B+1)T^3+(4A^2C^2+2AB^2C+4ABC\\
&+2AC-B^2)T^2+4AC^2(B+1)T+2AC^3)/B^2,\\
&4A^2C^2(AT^2+(B+1)T+C)((2A^3C-A^2B^2)T^4+(4A^2BC-2AB^3\\
&+4A^2C)T^3+(4A^2C^2-B^4+4ABC+2AC-B^2)T^2+(4ABC^2\\
&-2B^3C+4AC^2)T+2AC^3-B^2C^2)/B^3\bigg).\end{split}\] Put $T=1,$ and
note that the $X$-th coordinate of specialization at $T = 1$ of the
point $[2]P_1$ has the form
\[\frac{2AC((2AC-1)B^2+4AC(A+C+1)B+2AC(A+C+1)^2)}{B^2}.\]
A quick computation reveals that the remainder of the division of
the numerator by the denominator with respect to $B$ is equal
to\[4A^2C^2(A+C+1)(2B+A+C+1)\] and thus is non-zero provided
$AC\neq0$. By the Nagell-Lutz theorem (see
\cite[p.78]{Silverman-Tate}), $P_1$ is of infinite order on
$E_{1,1},$ the specialization of $E_1$ at $T = 1$, and thus $P$ is
of infinite order on $E_1$. Hence, the group $E_1(\mathbb{Q}(T))$ is
infinite.

Compute the points $[m]P$ on $E_1$ for $m=2,3,...$, next calculate
the corresponding point $\varphi_1^{-1}([m]P) = (v_m, z_m)$ on $C_1$
and solve the equation $f(u)=0$ for $u$. Put the calculated roots
into the expression for $u$, we get various $\mathbb{Q}(T)$-rational
solutions $(x,y,z)$ of Eq. (\ref{Eq11}) for $f(X)=AX+B+\frac{C}{X}$.

2) For $4AC-B^2=0,$ then $f(X)=\frac{(2AX+B)^2}{4AX},$ and
\[\begin{split}
\Delta(z)=&16A^2(2Az+B)^2\\
&\times(4A^2T^2z^2+(-8A^2BT^3-(8AB^2-4AB)T^2-2B^3T)z+B^2T^2).
\end{split}\]
If $4A^2T^2z^2+(-8A^2BT^3-(8AB^2-4AB)T^2-2B^3T)z+B^2T^2=w^2,$ then
$\Delta(z)$ is a square. This quadratic equation about $z,w$ can be
parameterized by
\[\begin{split}&z=\frac{2BT(4A^2T^2+2(2B-1)AT+B^2+t)}{4A^2T^2-t^2},\\
&w=\frac{BT(4(2t+1)A^2T^2+4(2B-1)AtT+2B^2t+t^2)}{4A^2T^2-t^2},\end{split}\]where
$t$ is a rational parameter. Then
\[\begin{split}u=&\frac{B(2AT-t)(4A^2T^2+(2B-1)2AT+B^2+t)}{2A(2AT+B)^2(2AT+t)},\\
or~&\frac{B(2AT+B)^2(2AT+t)}{2A(2AT-t)(4A^2T^2+2(2B-1)AT+B^2+t)}.\end{split}\]
So the rational parametric solutions of Eq. (\ref{Eq11}) are
\[\begin{split}
(x,y,z)=\bigg(&T,\frac{B(2AT-t)(4A^2T^2+2(2B-1)AT+B^2+t)}{2A(2AT+B)^2(2AT+t)},\\
&\frac{2BT(4A^2T^2+2(2B-1)AT+B^2+t)}{4A^2T^2-t^2}\bigg),\\
or~\bigg(&T,\frac{B(2AT+B)^2(2AT+t)}{2A(2AT-t)(4A^2T^2+2(2B-1)AT+B^2+t)},\\
&\frac{2BT(4A^2T^2+2(2B-1)AT+B^2+t)}{4A^2T^2-t^2}\bigg).\end{split}\]

Combining 1) and 2), we complete the proof of Theorem 1.1.
\end{proof}

\textbf{Example 1.} The point $[2]P$ on $E_1$ leads to the solutions
of Eq. (\ref{Eq11}):
\[\begin{split}
(x,y,z)=\bigg(&T,-\frac{(AT^2+BT+C)B}{(AT^2+(B+1)T+C)A},-\frac{(AT^2+(B+1)T+C)C}{BT}\bigg),\\
or~\bigg(&T,-\frac{(AT^2+(B+1)T+C)C}{(AT^2+BT+C)B},-\frac{(AT^2+(B+1)T+C)C}{BT}\bigg).
\end{split}\]

\begin{proof}[\textbf{Proof of Theorem 1.2.}]
For $f(X)=AX+\frac{C}{X}$, let \[x=T,y=u.\] Then Eq. (\ref{Eq11})
reduces to
\[(A^2T^2z+ACz)u^2+(-ATz^2-CT)u+ACT^2z+C^2z=0.\] If the above equation has
rational solutions $u$, then the discriminant
\[\Delta'(z)=A^2T^2z^4+(-4A^3CT^4-8A^2C^2T^2-4AC^3+2ACT^2)z^2+C^2T^2\]
should be a square.

Let us consider the curve $C_2:~v^2=\Delta'(z).$ If $AC\neq0,$ the
discriminant of $\Delta'(z)$ is non-zero as an element of
$\mathbb{Q}(T)$, then $C_2$ is smooth. By
the method described in \cite[p.77]{Mordell} (or see \cite[p.476,
Proposition 7.2.1]{Cohen1}), $C_2$ is birationally equivalent to
the elliptic curve
\[E_2:~Y^2=(X-2ACT^2)(X+2ACT^2)(X-2AC(2A^2T^4+(4AC-1)T^2+2C^2)).\]
Hence, if there exists a nonzero rational number $T$ such that the
elliptic curve $E_2$ has positive rank, then Eq. (\ref{Eq11}) has
infinitely many nontrivial rational solutions $(x,y,z)$.
\end{proof}

\textbf{Example 2.} When $A=C=1$, $T=4/3$, the elliptic curve $E_2$
leads to
\[E'_2:~U^2=V^3-2212V^2-82944V+183472128,\]
where $U=729Y,V=81X$. Then
\[z=\frac{U}{24(V-2212)},y=\frac{50(V-288)}{U}.\]
Note that the rank of $E'_{2}$ is one, so Eq. (\ref{Eq11}) has
infinitely many rational solutions $(4/3,y,z)$ for
$f(X)=X+\frac{1}{X}$, such as
\[\begin{split}(x,y,z)=&\bigg(\frac{4}{3},\frac{3}{5},\frac{9}{2}\bigg),\bigg(\frac{4}{3},\frac{231}{136},\frac{1309}{288}\bigg),
\bigg(\frac{4}{3},\frac{17849}{1945935},\frac{3076893}{13546}\bigg),\\
&\bigg(\frac{4}{3},\frac{102449900784}{62632357625},\frac{117165497209}{26287652160}\bigg).\end{split}\]

To get infinitely many positive rational solutions of Eq. (1.1), we need a theorem of Poincar\'e and Hurwitz (see
\cite[p.78]{Skolem}) about the density of rational points: If an elliptic
curve $E$ defined over $\mathbb{Q}$ has positive rank and at most
one torsion point of order two, then the set $E(\mathbb{Q})$ is dense
in $E(\mathbb{R})$. The same result holds if $E$ has three torsion
points of order two under the assumption that we have a rational
point of infinite order on the bounded branch of the set
$E(\mathbb{R})$.

If there is a point on $E'_2$ satisfying the condition
\[z=\frac{U}{24(V-2212)}>0,y=\frac{50(V-288)}{U}>0,\] then there are
infinitely many rational points on $E'_2$ satisfying it. Because the point
\[(V,U)=(8712,702000)\] on $E'_{2}$ satisfies the above condition, so Eq. (\ref{Eq11}) has infinitely many positive rational solutions
$(4/3,y,z)$ for $f(X)=X+\frac{1}{X}$.

\begin{proof}[\textbf{Proof of Theorem 1.3.}]
For $f(X)=AX+B+\frac{C}{X}$, let
\[y=xT^2,z=xT.\] Then Eq.
(\ref{Eq12}) becomes
\[g(x)=(T-1)^2(ABT^2x^2+(ACT^2+2ACT+AC)x+BC)=0.\]
Consider the above equation as a quadratic equation of $x$, if it
has rational solutions, the discriminant
\[\Delta(T)=A^2C^2T^4+4A^2C^2T^3+(6A^2C^2-4AB^2C)T^2+4A^2C^2T+A^2C^2\]
should be a square. Let us consider the curve $C_3:~v^2=\Delta(T).$

1) If $4AC-B^2\neq0,$ the discriminant of $\Delta(T)$ is non-zero, then $C_3$ is smooth. By the method
described in \cite[p.77]{Mordell} (or see \cite[p.476, Proposition
7.2.1]{Cohen1}), $C_3$ is birationally equivalent with the elliptic
curve
\[E_3:Y^2=(X+2A^2C^2)(X^2+4AC(AC-B^2)X+4A^3C^3(AC+2B^2)).\]
Because the map $\varphi_3: C_3\rightarrow E_3$ is complicated, we
omit it.

Note that $E_3$ contains the point
\[Q=(2A^2C^2, 8A^3C^3).\]
By the group law, we have
\[[2]Q=(-A^2C^2-2AB^2C+B^4, (AC-B^2)(A^2C^2+4AB^2C-B^4)),\]
\[\begin{split}
[4]Q=\bigg(&-(7A^8C^8+88A^7B^2C^7-420A^6B^4C^6+24A^5B^6C^5+50A^4B^8C^4\\
&+8A^3B^{10}C^3-20A^2B^{12}C^2+8AB^{14}C-B^{16})/(4(AC-B^2)^2(A^2C^2\\
&+4AB^2C-B^4)^2),\\
&(A^4C^4-20A^3B^2C^3+6A^2B^4C^2-4AB^6C+B^8)(A^8C^8+80A^7B^2C^7\\
&-180A^6B^4C^6+656A^5B^6C^5-282A^4B^8C^4-80A^3B^{10}C^3+76A^2B^{12}C^2\\
&-16AB^{14}C+B^{16})/(8(AC-B^2)^3(A^2C^2+4AB^2C-B^4)^3)\bigg).\end{split}\]

A quick computation reveals that the remainder of the division of
the numerator by the denominator of the $X$-th coordinate of $[4]Q$
with respect to $B$ is equal
to\[16A^3C^3(3B^{10}-17ACB^8+14A^2C^2B^6+26A^3C^3B^4-9A^4C^4B^2+A^5C^5)\]
and thus is non-zero provided $AC\neq0$. By the Nagell-Lutz theorem
(see \cite[p.78]{Silverman-Tate}), $[4]Q$ is of infinite order on
$E_3$, then there are infinitely many rational points on $E_3$ for
$4AC-B^2\neq0.$

For $m=2,3,...,$ compute the points $[m]Q$ on $E_3$, next calculate
the corresponding point $\varphi_3^{-1}([m]Q) = (v_m, T_m)$ on $C_3$
and solve the equation $g(x)=0$ for $x$. Put the calculated roots
into the expression for $x$, we get various rational solutions
$(x,y,z)$ of Eq. (\ref{Eq12}) for $f(X)=AX+B+\frac{C}{X}$.

2) For $4AC-B^2=0,$ then $f(X)=\frac{(2AX+B)^2}{4AX},$ and
\[\Delta(T)=A^2B^2(T^2+6T+1)(T-1)^2.\]
If $T^2+6T+1=w^2,$ then $\Delta(T)$ is a square. This quadratic
equation about $T,w$ can be parameterized by
\[T=-\frac{2(t-3)}{t^2-1},w=-\frac{t^2-6t+1}{t^2-1},\]where
$t$ is a rational parameter. Then
\[x=-\frac{(t+1)^2B}{16A},~or~-\frac{(t-1)^2B}{(t-3)^2A}.\]
So the rational parametric solutions of Eq. (\ref{Eq12}) are
\[\begin{split}(x,y,z)=\bigg(&-\frac{(t+1)^2B}{16A},-\frac{(t-3)^2B}{4(t-1)^2A},\frac{(t-3)(t+1)B}{8(t-1)A}\bigg),\\
or~\bigg(&-\frac{(t-1)^2B}{(t-3)^2A},-\frac{4B}{(t+1)^2A},\frac{2(t-1)B}{(t+1)(t-3)A}\bigg).\end{split}\]

Combining 1) and 2), then the proof of Theorem 1.3 is completed.
\end{proof}

\textbf{Example 3.} The point $-Q=(2A^2C^2, -8A^3C^3)$ on $E_3$
leads to the solutions of Eq. (\ref{Eq12}):
\[\begin{split}
(x,y,z)=\bigg(&-\frac{(AC-B^2)^2}{4A^2BC},-\frac{C}{B},\frac{AC-B^2}{2AB}\bigg),\\
or~\bigg(&-\frac{B}{A},-\frac{4ABC^2}{(AC-B^2)^2},\frac{2BC}{AC-B^2}\bigg);
\end{split}\]
the point $-[2]Q=(-A^2C^2-2AB^2C+B^4, -(AC-B^2)(A^2C^2+4AB^2C-B^4))$
on $E_3$ leads to the solutions of Eq. (\ref{Eq12}):
\[\begin{split}
(x,y,z)=\bigg(&-\frac{4ABC^2}{(AC-B^2)^2},-\frac{B(3AC-B^2)^2(AC+B^2)^2}{A(A^2C^2-6AB^2C+B^4)^2},\\
&\frac{2BC(3AC-B^2)(AC+B^2)}{(AC-B^2)(A^2C^2-6AB^2C+B^4)}\bigg),\\
or~\bigg(&-\frac{C(A^2C^2-6AB^2C+B^4)^2}{B(3AC-B^2)^2(AC+B^2)^2},-\frac{(AC-B^2)^2}{4A^2BC},\\
&\frac{(AC-B^2)(A^2C^2-6AB^2C+B^4)}{2AB(AC+B^2)(3AC-B^2)}\bigg).
\end{split}\]

\begin{proof}[\textbf{Proof of Theorem 1.4.}]
For $f(X)=AX+\frac{C}{X}$, Eq. (\ref{Eq12}) reduces to
\[(A^2y^2z^2+ACz^2)x^2+(-A^2yz^4-2ACyz^2-C^2y)x+ACy^2z^2+C^2z^2=0.\]
If the above equation has rational solutions $x$, then the
discriminant
\[\Delta(y)=-4A^3Cz^4y^4+(A^4z^8+4A^3Cz^6-2A^2C^2z^4+4AC^3z^2+C^4)y^2-4AC^3z^4\]
should be a square. Let us consider the curve
\[C_4:~v^2=-4A^3Cz^4y^4+(A^4z^8+4A^3Cz^6-2A^2C^2z^4+4AC^3z^2+C^4)y^2-4AC^3z^4.\]
Therefore, if there exist nonzero integers $A,C$ and a nonzero
rational number $z$ such that the elliptic curve $C_4$ has
infinitely many rational points $(y,v)$, then Eq. (\ref{Eq12}) has
infinitely many nontrivial rational solutions $(x,y,z)$.
\end{proof}

\textbf{Example 4.} When $A=1,C=-1$, $z=4$, $C_4$ is birationally
equivalent to \[E_4:~Y^2=X^3+48577X^2-4194304X-203746705408.\]Then
\[y=\frac{Y}{64(48577+X)},x=\frac{(48577+X)(X^2+97154X-450Y+4194304)}{(64X+3108928+Y)(64X+3108928-Y)}.\]
Note that the rank of $E_{4}$ is one, so Eq. (\ref{Eq12}) has
infinitely many rational solutions $(x,y,4)$ for
$f(X)=X-\frac{1}{X}$, such as
\[\begin{split}(x,y,z)=&\bigg(\frac{18}{13},\frac{1}{14},4\bigg),\bigg(\frac{18971}{1024},\frac{5891}{4061},4\bigg),
\bigg(\frac{825}{2204},\frac{1631}{10244},4\bigg),\\
&\bigg(\frac{58502431053824}{507945397025551},\frac{23133182812831}{48912132263775},4\bigg).\end{split}\]

In virtue of the theorem of Poincar\'e and Hurwitz (see
\cite[p.78]{Skolem}), $E_4$ has infinitely many rational points in
every neighborhood of any one of them. If there is a point on $E_4$
satisfying the condition
\[y=\frac{Y}{64(48577+X)}>0,x=\frac{(48577+X)(X^2+97154X-450Y+4194304)}{(64X+3108928+Y)(64X+3108928-Y)}>0,\] then there are
infinitely many rational points on $E_4$ satisfying it. Because the point \[(X,Y)=\bigg(\frac{112352}{49}, \frac{79764000}{343}\bigg)\] on $E_{4}$ satisfies the above condition, then
Eq. (\ref{Eq12}) has infinitely many positive rational solutions
$(x,y,4)$ for $f(X)=X-\frac{1}{X}$.

\section{Some related questions}

We have studied the rational parametric solutions of Eqs.
(\ref{Eq11}) and (\ref{Eq12}) for $f(X)=AX+B+\frac{C}{X}$, but we
don't get the same results for other simple Laurent polynomials.

\begin{question}
For $f(X)=AX^2+BX+C+\frac{D}{X}$ or $AX+B+\frac{C}{X}+\frac{D}{X^2}$, whether Eqs. (\ref{Eq11}) and (\ref{Eq12}) have rational solutions? If they have, are there infinitely many?
\end{question}

To find the integer solutions of Eqs. (\ref{Eq11}) and (\ref{Eq12})
is also an interesting question. Here we list some nontrivial
integer solutions of Eq. (\ref{Eq11}) for $f(X)=X+1+\frac{C}{X}$
with $1\leq C \leq 10$, $-100\leq x < y\leq 100$ and $-100\leq z\leq
100$ in Table 1.

\[\begin{tabular}{c|c}
\hline

 $C$ & $(x,y,z)$ \\

\hline

 $1$ & $(-3,-2,2)$  \\
\hline

 $2$ & $(-2, -1, 1),(-2, -1, 2)$  \\
 \hline

 $3$ & $(-4, -2, 8),(-4, 21, -84),(-3, -2, 6),(-2, 2, -12)$  \\
 \hline

 $4$ & $(-9, -3, 27),(-4, -3, 12),(-3, -1, 12)$  \\

 \hline

 $5$ & $(-16, -4, 64),(-5, -4, 20),(-4, -1, 20),(-2, 2, -20),(2, 3, 30)$  \\

\hline

 $6$ & $(-7, 2, -42),(-7, 3, -42),(-6, -5, 30),(-5, -1, 30)$  \\
 \hline

 $7$ & $(-7, -6, 42),(-6, -1, 42),(-5, 3, -35),(3, 13, 91)$  \\
 \hline

 $8$ & $(-8, -7, 56),(-7, -1, 56),(-5, 2, -40), (-5, 4, -40),(2, 11, 88),(4, 11, 88)$  \\
 \hline

 $9$ & $(-9, -8, 72),(-8, -1, 72)$  \\
 \hline

 $10$ & $(-10, -9, 90),(-9, -1, 90)$  \\
 \hline
\end{tabular}
\]
\begin{center}Table~~1. Some integer solutions of Eq. (\ref{Eq11})
for $f(X)=X+1+\frac{C}{X}$\end{center}

By some calculations, we can find a lot of Laurent polynomials such
that Eqs. (\ref{Eq11}) and (\ref{Eq12}) have integer solutions, but
it seems difficult to prove they have infinitely many integer
solutions. Hence, we raise

\begin{question}
Does there exist a Laurent polynomial $f(X)=AX+B+\frac{C}{X}$ with
$ABC\neq0$ such that Eqs. (\ref{Eq11}) and (\ref{Eq12}) have
infinitely many integer solutions?
\end{question}

\vskip20pt
\bibliographystyle{amsplain}

\end{document}